\documentclass[12pt]{article}
\usepackage{amsmath, amsthm}
\usepackage{amsfonts,amssymb}
\usepackage{epsfig}
\usepackage{graphicx}
\usepackage{subfigure}
\newcommand{\R}{{\mathbb R}}

\newcommand{\mb}{\mathbf}
\newcommand{\B}{\bigskip}
\newcommand{\m}{\medskip}
\newtheorem{theorem}{Theorem}[section]

\newtheorem{lemma}[theorem]{Lemma}
\newtheorem{prop}[theorem]{Proposition}

\theoremstyle{definition}
\newtheorem{definition}[theorem]{Definition}
\newtheorem{example}[theorem]{Example}

\numberwithin{equation}{section}
\def\sqr#1#2{{\vcenter{\vbox{\hrule height.#2pt
 \hbox{\vrule width.#2pt height#1pt \kern#1pt
 \vrule width.#2pt}
 \hrule height.#2pt}}}}
\def\square{\mathchoice\sqr68\sqr68\sqr{2.1}3\sqr{1.5}3}
\date{ }

\title{The Entropy of an Overlapping Dynamical System} 
\author{Michael  Barnsley \\  The
Australian National University \\Canberra, Australia \\ \tt michael.barnsley@anu.edu.au
\\ \\
Brendan Harding \\  The
Australian National University \\Canberra, Australia  \\ \tt brendan.harding@anu.edu.au
\\ \\
Andrew Vince \\ Department of Mathematics, University
of Florida \\ Gainesville, FL, USA \\ \tt avince@ufl.edu }

\begin{document}
\maketitle

\begin{abstract}  
The term ``overlapping" refers to a certain fairly simple type of piecewise continuous function from the unit interval to itself and also to a fairly simple type of iterated function system (IFS) on the unit interval.  A correspondence between these two classes of objects is used (1) to find a necessary and sufficient condition for a fractal transformation from the attractor of one overlapping  IFS to the attractor of another overlapping IFS to be a homeomorphism and (2) to find a formula for the topological entropy of the dynamical system associated with an overlapping function. 
\end{abstract}


\section{Introduction}  \label{secIntro}
 Iterated maps on an interval provide the simplest examples of dynamical systems. Parameterized families of geometrically simple continuous dynamical systems on an interval have a rich history because of their intricate behaviour, the insights they provide into higher dimensional systems, and diverse applications. Numerous papers have been written concerning their invariant measures, entropies, and behaviours; we note in particular the works of Collet and Eckman \cite{collet}, and Milnor and Thurston \cite{ milnor}.  Many results in the literature concern the topological entropy of continuous systems \cite{adler,mis}, but piecewise continuous maps have also received attention, \cite{hoffbauer,parry,renyi}.  

\begin{figure}[htb] \label{fA}
\begin{center} 
\includegraphics[width=5.5486in, keepaspectratio]{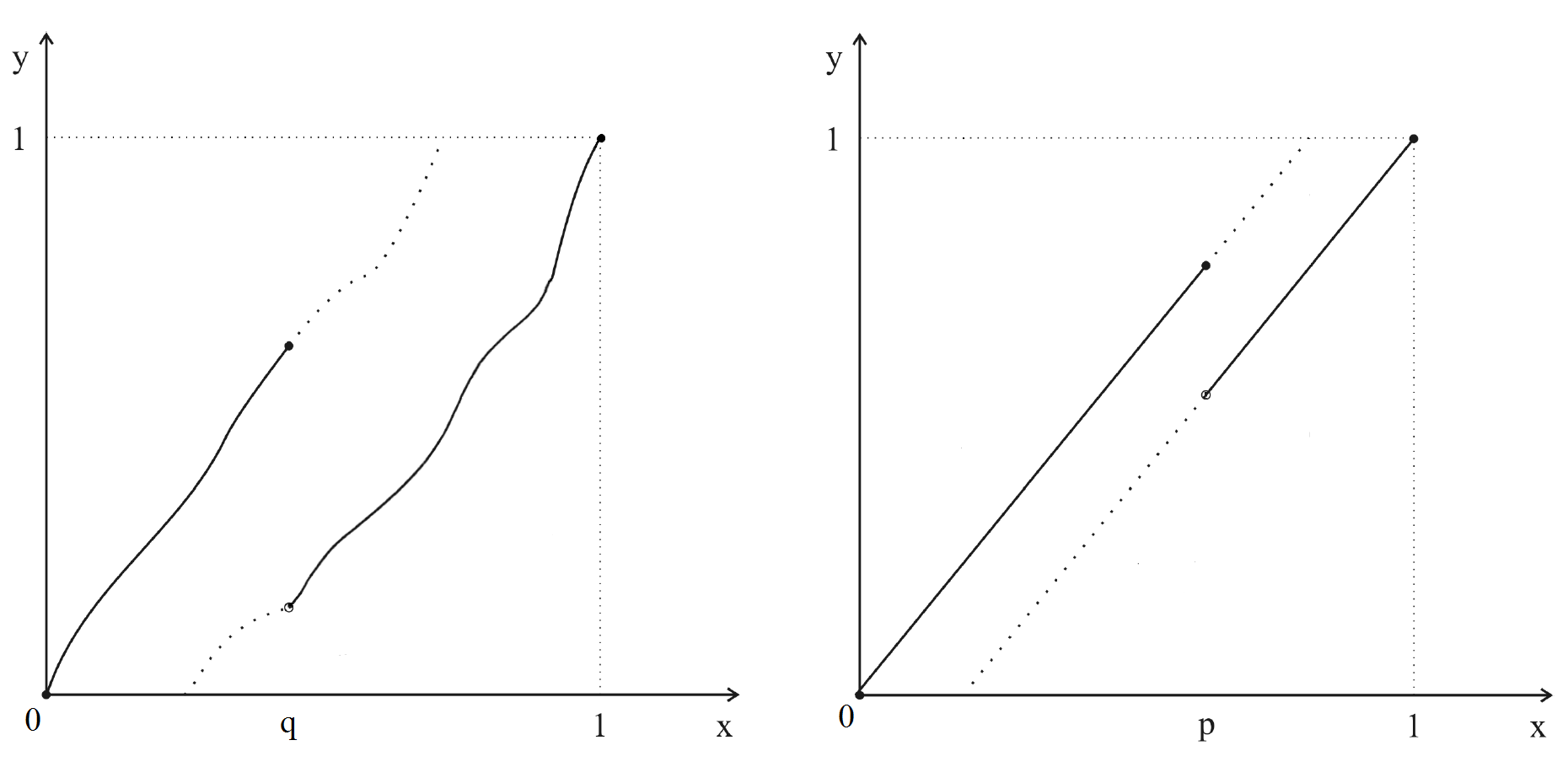}
\caption{An overlapping dynamical system (left) and the unform case (right).}
\end{center}
\end{figure}

 One of the main results of this paper is a formula for the topological entropy of a dynamical system $([0,1], T)$, where $T$ is a piecewise continuous function from the interval $[0,1]$ onto itself consisting of two continuous pieces, as shown on the left in Figure~\ref{fA}.   More precisely we are interested in functions $T \, : \, [0,1]\rightarrow [0,1]$  of the form
\begin{equation} \label{eqT1}
T(x) = \begin{cases} g_0(x) \quad  &\text{if $0\leq x < q$}  \\  
g_1(x) \quad &\text{if $q \leq x \leq 1$}, \end{cases}\end{equation}
or
\begin{equation} \label{eqT2}
T(x) = \begin{cases} g_0(x) \quad  &\text{if $0\leq x \leq q$}  \\  
g_1(x) \quad &\text{if $q < x \leq 1$}, \end{cases}\end{equation}
where $g_0 \, : \, [0,1] \rightarrow [0,1]$ and $g_1\, : \,[0,1] \rightarrow [0,1]$ are continuous increasing functions such that
\begin{enumerate}
\item $g_0(0) = 0,  \, g_1(1) = 1$, 
\item $0 < g_1(q) < g_0(q) < 1$ for some $q\in (0,1)$, and
\item $g_0$ and $g_1$ are expansive, i.e. there is an $s > 1$ such that $|g_i(x)- g_i(y)| \geq s|x-y|$ for $i=0,1$ and for all $x\in [0,1]$. 
\end{enumerate}
Call such a dynamical system an {\bf overlapping} dynamical system.   It is ``overlapping" in the sense
that $g_0([0,q))\cap g_1((q,1]) \neq \emptyset$.
Note that the expansive condition 3 holds, for example, if  $g_0,g_1$ are differentiable and there is an  $s > 1$ such 
that $g_0'(x) \geq s$ and $g_1'(x) \geq s$ for all $x\in [0,1]$. \m

Associated with the dynamical system $T$ there are two special itineraries, called {\it critical itineraries}
$\alpha := \alpha_0, \alpha_1, \alpha_2 \dots$ and $\beta := \beta_0, \beta_1, \beta_2, \dots$, 
where $\alpha_n, \beta_n \in \{0,1\}$ for all $n\geq 0$ (see Definition~\ref{defMaskSeq}). Theorem~\ref{thmMain}, states that the topological entropy of $T$ is
$-\ln r$ where $r$ is the smallest solution $x\in (0,1)$ to the equation
$$\sum_{n=0}^{\infty} \alpha_n x^n = \sum_{n=0}^{\infty} \beta_n x^n.$$ 

The proof of this theorem relies on finding a ``uniform"  dynamical system that is topologically conjugate to 
the dynamical system $([0,1],T)$.  By uniform we mean a function $U$  of the form shown on the right in Figure~\ref{fA}, where the two branches are lines of equal slope.  For such a dynamical system it is well known that the entropy is $\ln r$, where $r$ is the slope of the lines.   That there exists such a topologically conjugate uniform dynamical system follows from \cite[Theorem 1]{denker} and can also be deduced from \cite{parry}. 
What was not known prior to this work, is the explicit relationship between $T$, on the left in Figure \ref{fA}, and the parameters $p$ and $r$ that uniquely determine $U$, on the right in Figure \ref{fA}. In this paper
(Theorem~\ref{thmMain}),  we construct such a topologically conjugate $U$ by determining the parameters $p$ and $r$ in terms of just the two critical itineraries $\alpha$ and $\beta$ of $T$.  Our approach is constructive in character.  We make use of an analogue of the kneading determinant of \cite{milnor}, appropriate for discontinuous interval maps, and thereby avoid a measure-theoretic existential proof such as those in \cite{hoffbauer, parry}. \m

Although some main results  concern dynamical systems, the underlying subject of the paper is a surprising connection between two areas -  the dynamics of a single map, on the one hand, and iterated function systems, on the other.  Our proof of the main theorems in this paper depends on this correspondence.  The correspondence  is such that two dynamical systems are topologically conjugate if and only if the attractors of  the two corresponding iterated function systems are related by a fractal homeomorphism. 
Indeed, one motivation for undertaking this research was our desire to establish, and to be able to compute, fractal homeomorphisms between attractors of  iterated function systems - for applications such as those in \cite{BHI}. \B

\begin{figure}[htb]  \label{fB}
\begin{center}
\includegraphics[width=4.5cm, keepaspectratio]{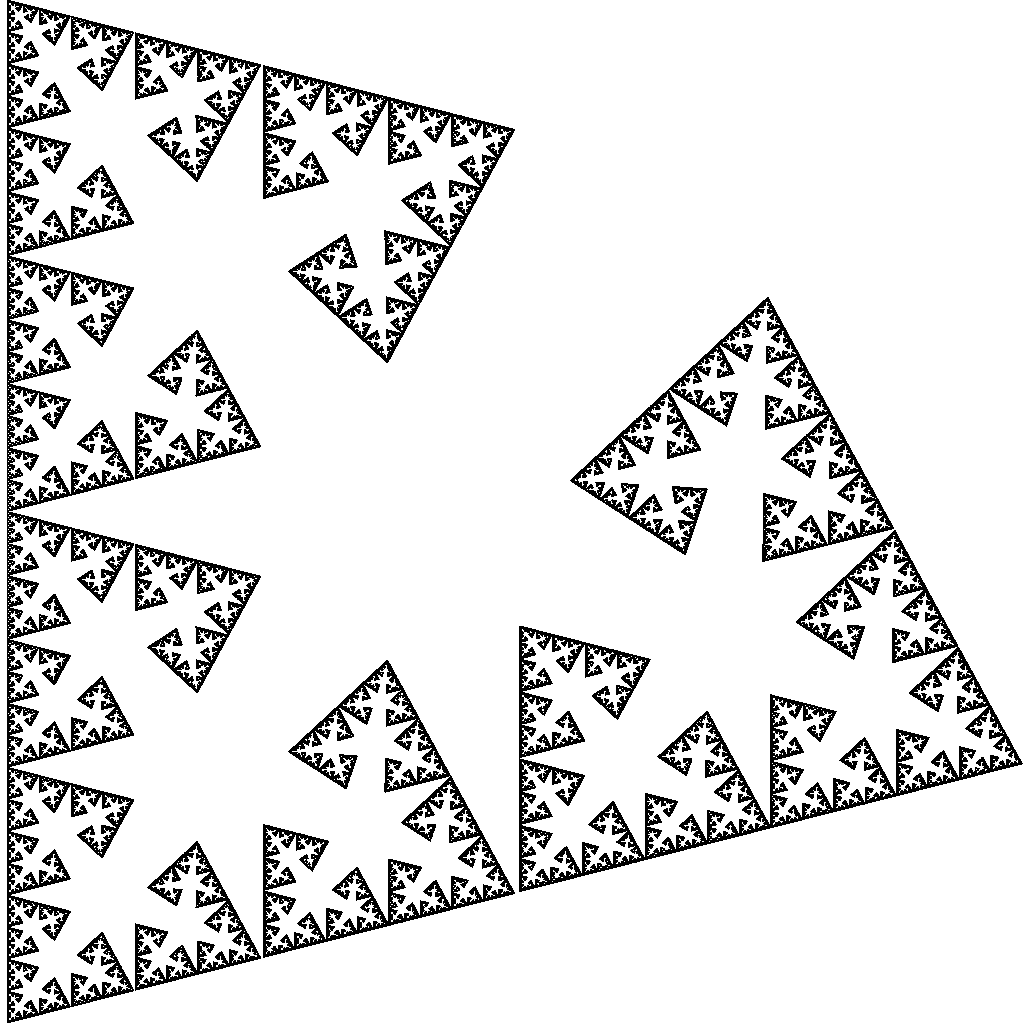}  \hskip 3cm 
\includegraphics[width=4.5cm, keepaspectratio]{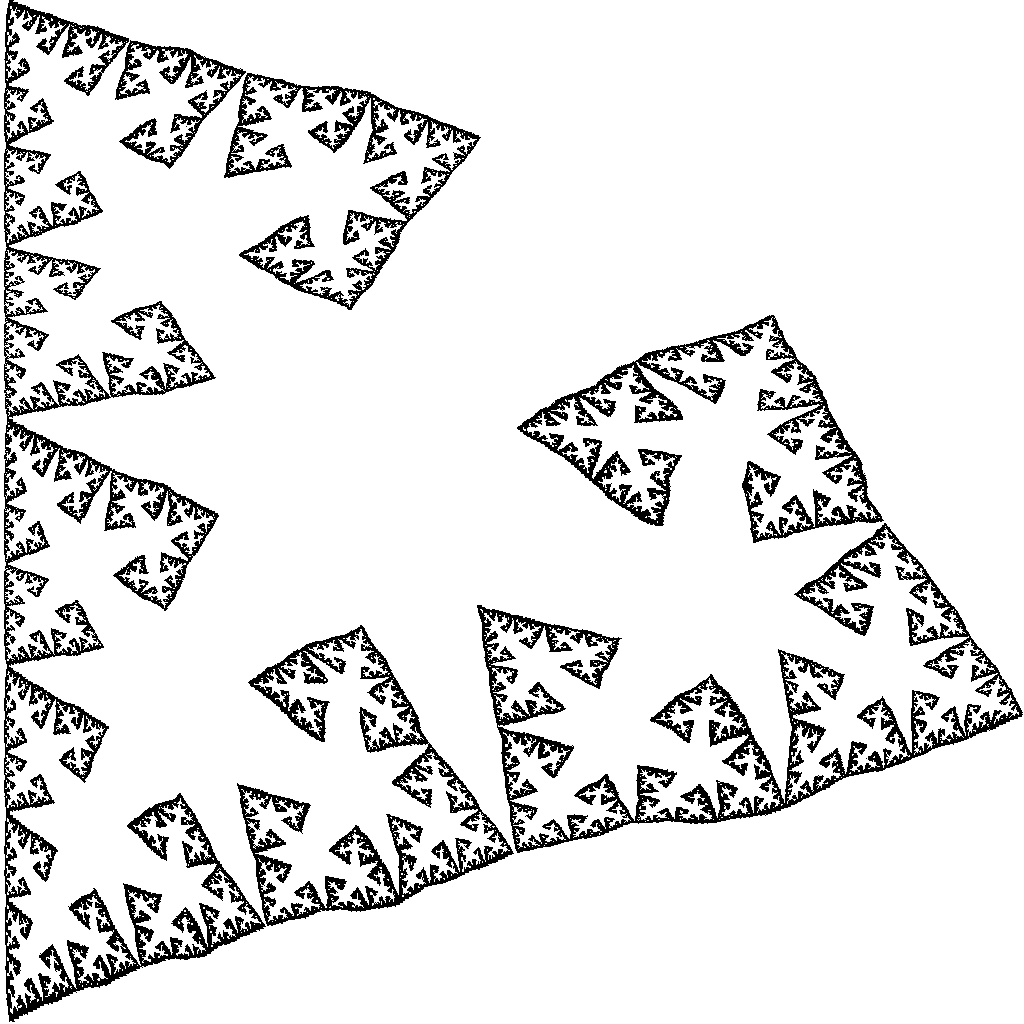} 
\vskip 3mm
\caption{The attractor of an IFS (left) and its image under a fractal homeomorphism (right).}
\end{center} 
\end{figure}

An iterated function system (IFS) is a standard method for constructing a self-referential fractal, the attractor of the IFS   usually being a fractal.  Given two iterated function systems  with the same number of functions, a method for transforming the attractor of one to the attractor of the other has been laid out in \cite{monthly}.  Figure 2 shows such a fractal transformation.  Even if the attractors themselves are mundane, the fractal transformations between them may be interesting.  In Figure 3, for example, the attractors are simply the unit square $\square$.  To visualize the fractal
transformation we can observe its effect on a ``picture".    By {\it picture} we mean a function $c \, : \, \square \rightarrow {\cal C}$, where $\cal C$ denotes the color palate, for example
 ${\cal C} = \{ 0,1,2,\dots, 255\}^3$.  A fractal transformation $h \, : \, \square \rightarrow \square$  induces a map from a picture on one attractor to a picture on the other attractor given  by
 $$h(c) := c \circ h.$$  The particular fractal transformation depicted in Figure 3 is a homeomorphism.  The question of when a fractal transformation is a homeomorphism, difficult in even simple situations, is answered in this
paper for the case of an ``overlapping" IFS on the unit interval, i.e. for an IFS $([0,1]; \, f_0,f_1)$  consisting of two contractions $f_0,f_1$ defined on the unit interval $[0,1]$ such that $[0,1] = f_0([0,1] \cup f_1([0,1]$  and $f_0([0,1] \cap f_1([0,1] \neq \emptyset$.   One necessary and sufficient condition is proved as part of Theorem~\ref{thmFT-TC}:  a fractal transformation is a homeomorphism if and only if the critical itineraries $\alpha$ and $\beta$ associated with one IFS equal the critical itineraries associated with the other.  \B

\begin{figure}[htb] \label{fC}
\begin{center}
\includegraphics[width=5cm, keepaspectratio]{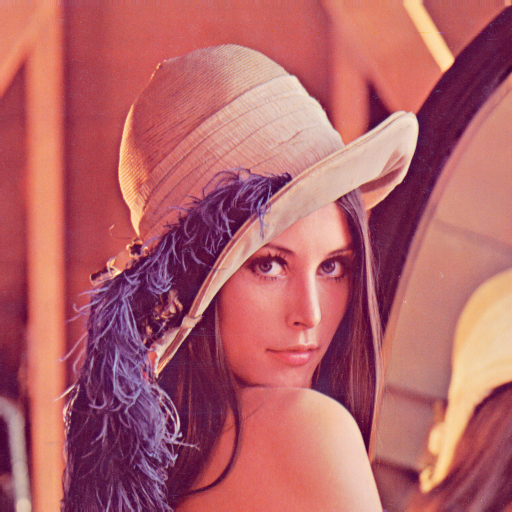} \hskip 3cm 
\includegraphics[width=5cm, keepaspectratio]{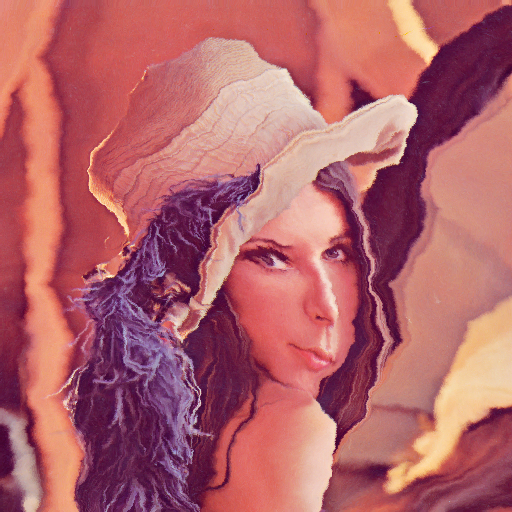}
\vskip 2mm
\caption{A fractal homeomorphism applied to the original picture.}
\end{center} 
\end{figure}

The organization of the paper is as follows.  Basic definitions and facts about iterated function systems and their attractors are reviewed in Section~\ref{secADS}.  The dynamical system associated with an IFS is also defined 
in that section.  Fractal transformations and how they are constructed using masks (Theorems~\ref{thmMask} and \ref{thmPi-Tau}) are the subjects of Section~\ref{secFT}. The particular type of IFS that is central to this paper, 
an overlapping IFS, is defined in Section~\ref{secOverlap}.   A uniform IFS, a special case of an overlapping IFS, is also discussed in that section.  Each point of the attractor of an IFS can be assigned an address.   The the address space of the attractor of an overlapping IFS is the topic of  Section~\ref{secAddSpace}.  The two critical itineraries are defined in this section, and two important results
are stated.  Theorem~\ref{thmAddSpace} characterizes the address space of an overlapping IFS in terms of the critical itineraries.  Theorem~\ref{thmFT-TC} states that the following four conditions are equivalent: (1) the address spaces of two overlapping IFSs are equal; (2) the corresponding critical itineraries are equal; (3)  the two IFS are related by a fractal homeomorphism; and (4)  the two associated dynamical systems are topologically conjugate. Theorems~\ref{thmAddSpace} and \ref{thmFT-TC}  lead to the main result on topological entropy, Theorem~\ref{thmMain}, stated and proved, with the aid of several lemmas, in Section~\ref{secEntropy}.  

\section{An IFS and its Associated Dynamical System} \label{secADS}

Basic results on iterated function systems and their associated dynamical systems are contained in this section. 
We begin in the setting of a complete metric space and specialize to the unit interval on the real line in
Section~\ref{secOverlap}.  

Let $\mathbb{X}$ be a complete metric space. If $f_{m}:\mathbb{X}
\rightarrow\mathbb{X}$, $m=1,2,\dots,M,$ are continuous maps, then
$\mathcal{F}=\left(  \mathbb{X};f_{1},f_{2},...,f_{M}\right)$ is called an
{\it iterated function system} (IFS). To define the attractor of an IFS, 
first define
\[
\mathcal{F}(B)=\bigcup_{f\in\mathcal{F}}f(B)
\]
for any $B\subset\mathbb{X}$. By slight abuse of terminology we use the same
symbol $\mathcal{F}$ for the IFS, the set of functions in the IFS, and for the
above map. For $B\subset\mathbb{X}$, let $\mathcal{F}^{k}(B)$ denote the
$k$-fold composition of $\mathcal{F}$, the union of $f_{i_{1}}\circ f_{i_{2}%
}\circ\cdots\circ f_{i_{k}}(B)$ over all finite words $i_{1}i_{2}\cdots i_{k}$
of length $k.$ Define $\mathcal{F}^{0}(B)=B.$ A nonempty compact set 
$A\subset\mathbb{X}$ is said to be an {\it attractor} of the IFS $\mathcal{F}$ if
\begin{enumerate}
\item $\mathcal{F}(A)=A$ and
\item $\lim_{k\rightarrow\infty}\mathcal{F}^{k}(B)=A,$ for all compact sets
$B\subset \mathbb X$, where the limit is with respect to the Hausdorff metric.
\end{enumerate}
\m

A function $f:\mathbb{X}\rightarrow\mathbb{X}$ is called a
{\it contraction} with respect to a metric $d$ if there is an $0\leq s <1$
such that $d(f(x),f(y))\leq s \,d(x,y)$ for all $x,y\in{\mathbb{R}}^{n}$.  An IFS with the property that each function is a contraction will be called a {\it contractive} IFS.  In his seminal paper  Hutchinson \cite{hutchinson} proved  that 
a contractive IFS on a complete metric space has a unique attractor. 

For a contractive IFS, it is possible to assign to each point of the attractor an ``address" as follows.  
 Let  $\Omega = \{1,2, \dots ,N\}^{\infty}$ denote the set of infinite strings using symbols $1,2, \dots , N$.  For a
string $\omega \in \Omega$, denote the $n^{th}$ element, $n \geq 0$,  in the string by $\omega_n$, and 
denote by $\omega|_n$ the string consisting of the first $n+1$ symbols in $\omega$, i.e., $\omega|_n  = \omega_0\omega_1\cdots \omega_{n}$.   Moreover,
 we use the notation $$f_{\omega|_n} := f_{\omega_0} \circ f_{\omega_1}\circ \cdots \circ f_{\omega_n}.$$
The set $\Omega$ can be given the product topology induced from the discrete topology on $\{1,2, \dots , N\}$. The product topology  is the same as the topology induced by the metric $d(\omega
,\sigma)=2^{-k}$ where $k$ is the least index such that $\omega_{k}\neq \sigma_{k}$.  The space $(\Omega,d)$ is a compact metric space.

\begin{definition} \label{defCode}
Let $\mathcal{F}=\left(  \mathbb{X};f_{1},f_{2},...,f_{N}\right)$  be a contractive IFSs on a complete metric space $\mathbb{X}$ with attractor $A$.  
The map $\pi: \Omega \rightarrow A$ defined by
 $$\pi(\sigma) :=
\lim_{k\rightarrow\infty}f_{\sigma|_k}(x)$$
is called the {\bf coding map} of $\mathcal F$.  
\end{definition}

For a contractive IFS  it is well known \cite{hutchinson} that the limit exists and is independent of $x \in \mathbb X$.  Moreover $\pi$  is continuous, onto, and satisfies the following commuting diagram  for each $n = 1,2,\dots, N$.
\begin{equation*} \label{diagramCode}
\begin{array}
[c]{ccc}
\Omega & \overset{s_{n}}{\rightarrow} & \Omega\\
\pi\downarrow\text{\ \ \ \ } &  & \text{ \ \ \ }\downarrow\pi\\
\mathbb X & \underset{f_{n}}{\rightarrow} & \mathbb X
\end{array}
\end{equation*}
The symbol $s_{n}:\Omega\rightarrow\Omega$ denotes the inverse shift map defined by $s_{n}(\sigma)=n \sigma.$

\begin{definition}\label{defSection}
 A {\bf section} of the coding map $\pi$ is a function $\tau \, : \, \Omega \rightarrow A$ 
such that $\pi \circ \tau$ is the identity.  For $x \in A$, the string $\tau (x)$  is referred to as the {\bf address} of $x$.  Call the set $\Omega_{\tau} := \tau(A)$ the {\bf address space} of the section $\tau$.   
\end{definition}

\begin{definition}\label{defShift}
Let $S$ denote the shift operator on $\Omega$, i.e, $S(n \sigma) = \sigma$ for any $n \in \{ 1,2, \dots , N\}$ and any $\sigma \in \Omega$.  
A subset $W \subseteq \Omega$ will be called {\bf shift invariant} if $S(W) \subseteq W$.  
If $\Omega_{\tau}$ is shift invariant, then $\tau$ is called a {\bf shift invariant section}. 
\end{definition}

 The following example demonstrates the naturalness of shift invariance. 

\begin{example} \label{EX1}  Consider the IFS ${\mathcal F} = ( \R \, ; \, f_0,f_1)$ where $f_0(x) = \frac12\, x$ and $f_1(x) = \frac12\, x + \frac12$.  The attractor is the interval $[0,1]$.  An address of a point $x$ is a binary representation of $x$.  In choosing a section $\tau$ one must decide, for example, whether to take $\tau (\frac14) = .01$ or $\tau (\frac14) = .00111\cdots$. If the section $\tau$ is shift invariant, this would imply, for example, that if $\tau(\frac14) = .00111\cdots$, then  $\tau (\frac12) = .0111\cdots$, not  $\tau (\frac12) = .100\cdots$.
\end{example} 

Call an IFS {\it injective} if each function in the IFS is injective.  Theorem~\ref{thmMask} below, which is proved in \cite{BV}, states that every shift invariant section of an injective IFS can be obtained from a mask.

\begin{definition} \label{defMask}
For an IFS $\mathcal F$ with attractor $A$, a {\bf mask} is a partition $M = \{M_i, 1\leq i \leq N\}$ of $A$ such that $M_i \subseteq f_i(A)$ for all $f_i\in {\mathcal F}$.   
\end{definition}

\begin{definition}  \label{defSM}  Given an injective IFS $\mathcal F$ with a mask $M = \{M_i, 1\leq i \leq N\}$, the {\bf section} $\mb{\tau_M}$ {\bf associated with mask} $\mb M$ is the function $\tau_M \, : \, A \rightarrow \Omega$ defined as follows.  
Let $\Omega_k$ denote the set of all finite strings of length $k$ in the symbols $\{1,2,\dots , N\}$.  
For each $k\geq 0$ define a partition $M^k = \{ M_{\sigma} \, : \, \sigma \in \Omega_k\}$ of $A$ recursively by taking $M^1 = M$ and $$M^{k+1} =\{ M_{\sigma \, j} = M_{\sigma} \cap f_{\sigma}(M_j) \, : \, \sigma \in \{1,2, \dots, N\}^k, 1\leq j \leq N \}.$$  A straightforward induction shows that $M^k$ is indeed a partition of $A$ for every $k\geq 0$, and that each such partition is a refinement of the previous partition, in particular $M_{\sigma \, j} \subseteq M_{\sigma}$ for all finite $\sigma$ and all
$j \in \{1,2, \dots, N\}$.    (Note that, for some values of $\sigma$, the sets $M_{\sigma}$ may be empty.)  Moreover, since the functions in the IFS are contractions, the maximum diameter of the sets in $M^k$ approachs $0$ as $k \rightarrow \infty$.  Since each $x \in \mathbb X$ lies in a unique nested sequence $$M_{i_0} \supseteq   M_{i_0\, i_1 } \supseteq  M_{i_0\, i_1 \, i_2 } \supseteq \cdots,$$
we can define $\tau_M (x) = i_0\, i_1 \, i_2 \cdots$. Note that this definition of $\tau_M$ is equivalent to saying that
$$x \in f_{i_0} \circ f_{i_1} \circ f_{i_2}\circ \cdots \circ f_{i_{k-1}} (M_{i_k})$$
for all $k\geq 0$.  That $\tau_M$  is indeed a section is part of Theorem~\ref{thmMask} below, whose proof 
appears in \cite{BV}.
\end{definition}

\begin{lemma} \label{lemM} With notation as above, for any injective IFS and for any finite string $\sigma$ and symbol $j$, 
we have $M_{j\sigma} = M_j\cap f_j(M_{\sigma})$. 
\end{lemma}

\proof The result will be proved by induction on the length of $\sigma$. Concerning length $1$,  it is easy to check, from the definition of the partition, that $M_{ji} = M_j\cap f_j(M_{i})$. Now 
$$\begin{aligned} M_{j\sigma i} &= M_{j\sigma} \cap f_{j\sigma}(M_i) = M_j \cap f_j(M_{\sigma}) \cap f_{j\sigma}(M_i) \\ &=  M_j \cap f_j(M_{\sigma}) \cap f_j(f_{\sigma}(M_i)) =
M_j \cap f_j(M_{\sigma} \cap f_{\sigma}(M_i)) = M_j \cap f_j(M_{\sigma i}),\end{aligned}$$
the second to last
equality using that $f_j$ is injective. 
\qed

\begin{theorem} \label{thmMask} Let $\mathcal F$ be a contractive and injective IFS. 
\begin{enumerate}
\item If $M$ is a mask, then $\tau_M$ is a shift invariant section of $\pi$. 
\item If $\tau$ is a shift invariant section of $\pi$, then $\tau = \tau_M$ for some mask $M$.  
\end{enumerate} 
\end{theorem} 

\begin{definition} \label{defItinerary}
Let $\mathcal F$ be an injective IFS with attractor $A$.  Given a mask $M$ for $\mathcal F$, define a function $T_{({\mathcal F},M)} \, : \, A \rightarrow A$ by $$T_{({\mathcal F},M)}(x) := f_i^{-1} (x) \quad \text{when} \quad x\in M_i.$$  The pair $(A, T_{({\mathcal F},M)})$ will be called the {\bf dynamical system associated with} $\mathcal F$ {\bf and} $\mb M$.  The {\bf itinerary} of a point $x\in A$ is the string $i_0\, i_1 \, i_2 \cdots \in \Omega$, where $i_k$ is the unique integer, $1\leq i_k \leq N$, such that $$T_{({\mathcal F},M)}^k(x) \in M_{i_k}.$$  
\end{definition}

\begin{prop}  \label{propItinerary}  If $\mathcal F$ is an injective masked IFS with associated dynamical system \linebreak $(A, T_{({\mathcal F},M)})$, then, for all $x\in A$, the itinerary of $x$ is its address $\tau_M(x)$.
\end{prop}

\proof By its definition, $i_0\, i_1 \dots$ is the itinerary of $x$ if and only if $f_{i_{k-1}}^{-1}\circ \cdots \circ f_{i_1}^{-1}\circ f_{i_0}^{-1}(x) \in M_{i_k}$ for all $k \geq 0$.  But this 
is equivalent to $x \in f_{i_0}\circ f_{i_1} \circ \cdots \circ f_{k-1}(M_{i_k})$ for all $k \geq 0$, which, as noted above, defines the sections. \qed

\section{Fractal Transformation} \label{secFT}

Consider two contractive IFSs $\mathcal{F}=\left(  \mathbb{X};f_{1},f_{2},...,f_{N}\right)$ and $\mathcal{G}=\left(  \mathbb{Y};g_{1},g_{2},...,g_{N}\right)$  with the same number $N$ of functions on  complete metric spaces
$\mathbb{X}$ and $\mathbb{Y}$.  Basically a fractal transformation from  ${\mathcal F}$ to ${\mathcal G}$ is  a map $h \, : \, A_F \rightarrow A_G$ that sends a point in the attractor $A_F$ of $\mathcal F$ to the point in the attractor $A_G$ of $\mathcal G$ with the same address.  More specifically:

\begin{definition} \label{defFT}
 Let $A_F$ and $A_G$ be the  attractors and $\pi_F$ and $\pi_G$ the coding maps of 
contractive IFSs $\mathcal F$ and $\mathcal G$, respectively.  
A map $h \, : \, A_F \rightarrow A_G$ is called a {\bf fractal transformation } if there exist shift invariant sections $\tau_F$ and $\tau_G$ such that the following diagram commutes: 
\begin{equation} \label{diagFT}
\begin{array}
[c]{ccc}
A_F & \overset{h}{\rightarrow} & A_G \\ \text{$\tau_F$} \searrow & & \swarrow \text{$\tau_G$} \\
 & \Omega & 
\end{array}
\end{equation}
i.e., the transformation $h$ takes each point $x\in A_F$ with address $\sigma = \tau_F(x)$ to the point $y \in A_G$ with the same address $\sigma =\tau_G(y) $.  A fractal transformation that is a homeomorphism is called a 
{\bf fractal homeomorphism}.  
\end{definition}

Theorem~\ref{thmPi-Tau}, proved in \cite{BV},  states that the fractal transformations between $A_F$ and $A_G$ are exactly maps of the form $\pi_G\circ \tau_F$ or $\pi_F\circ\tau_G$ for some shift invariant sections $\tau_F, \, \tau_G$.  

\begin{theorem}  \label{thmPi-Tau} Let $\mathcal F$ and $\mathcal G$ be contractive IFSs.  With notation as above
\begin{enumerate} 
\item If $h \, : \, A_F \rightarrow A_G$ is a fractal transformation with corresponding sections $\tau_F$ and $\tau_G$, then  $h = \pi_G \circ \tau_F$ and $h^{-1} = \pi_F \circ \tau_G$ .
\item If $\tau_F$ is a shift invariant section for $\mathcal F$, then $h := \pi_G \circ \tau_F$ is a fractal transformation.
\end{enumerate}
\end{theorem}

\section{Overlapping IFS} \label{secOverlap}

The type of IFS that is the subject of this paper is what will be called an overlapping IFS on the unit interval of
the real line.

\begin{definition} \label{defOverlap}
An {\bf overlapping IFS} is an IFS $${\mathcal F} = ([0,1] ; f_{0}(x), f_{1}(x)), $$ where the functions are 
continuous, increasing, contractions that satisfy $$ f_0(0) = 0, \qquad f_1(1) = 1, \qquad 0 < f_1(0) < f_0(1) < 1.$$
\end{definition}

The attractor of an overlapping IFS  is the unit interval $[0,1]$.  If it were the case that $ f_1(0)  =  f_0(1)$, then the two sets $f_0([0,1]$ and $f_1([0,1])$ would be ``just touching"; since  $ f_1(0)  <  f_0(1)$,  they are ``overlapping". \m

Next we fix some notation used in the remainder of the paper.  
The coding map for $\mathcal F$ will be denoted by $\pi := \pi_{\mathcal F}$. We consider masks for an overlapping IFS $\mathcal F$ of the form $$M_q^+ = \{\, [0,q), [q, 1] \,\} \text{ or }  M_q^- = \{ \, [0,q], (q, 1]\, \}, \text{ where } f_1(0) < q  <f_0(1).$$   

\begin{definition}  \label{defMaskP}   
The point $q$ will be called the {\bf mask point}. 
\end{definition}

For a masked overlapping IFS let  $\tau^{+}_{q}$ and $\tau^{-}_{q}$ denote the sections corresponding to $M_q^+$ and $M_q^-$, respectively. The two  respective address spaces are denoted by 
$$\Omega_q^- = \tau_q^-( [0,1]),  \qquad \text{and} \qquad \Omega_q^+ = \tau_q^+( [0,1]).$$
For a masked overlapping IFS, the associated dynamical systems, as defined in the previous section,  are
 $([0,1],  T_q^+)$ and  $([0,1],  T_q^-)$, where
$$T_q^+(x) =  \begin{cases} f_0^{-1} \quad \text{if} \quad  x < q \\  f_1^{-1} \quad \text{if} \quad x \geq q \end{cases} \qquad \text{and} \qquad
 T_q^-(x) =  \begin{cases} f_0^{-1}  \quad \text{if} \quad  x \leq  q 
\\  f_1^{-1} \quad \text{if} \quad x > q.
 \end{cases}$$ 
Since $f_0$ and $f_1$ are contractions,
the inverses $g_0 = f_0^{-1}$ and $g_1 = f_1^{-1}$ are expansions.  Also since $f_1(0) < q  <f_0(1)$, we have $0 < g_1(q) < g_0(q) < 1$ .  Therefore the dynamical system associated with an overlapping IFS is
an overlapping dynamical system as defined and discussed in the introduction. We will refer to such a dynamical system as an {\bf overlapping dynamical system}.  So there is a bijection between overlapping dynamical systems and
masked overlapping iterated function systems. \m

For our purposes, the following can serve as a definition of the entropy of a dynamical system.   Note that 
$|\Omega_{q,n}^{+}|= |\Omega_{q,n}^{-}|$, so this definition is consistent with the one in \cite{parry}.

\begin{definition} \label{defEntropy} The {\bf topological entropy} $h( T_{q}^{\pm})$ of an overlapping dynamical system $([0,1],  T_{q}^{\pm})$ is 
\[
h(T_{q}) = \lim_{n \to \infty} \frac{1}{n} \log |\Omega_{q,n}^{+} | =  \lim_{n \to \infty} \frac{1}{n} \log |\Omega_{q,n}^{-} | ,
\]
where $\Omega_{q,n}^{\pm} := \{ \omega|_n \,:\, \omega \in \Omega^{\pm}_{q} \}$.  
\end{definition}
 
The following  special case of an overlapping IFS plays an important rolel.  

\begin{definition}  The IFS ${\mathcal U}_{a} = ([0,1] ; L_{0}(x), L_{1}(x))$ where
\begin{equation} \label{eqAffine} \begin{aligned} L_0 (x) &= ax. \\  L_1 (x) &= ax + 1-a 
\end{aligned} \end{equation} 
will be called a {\bf uniform} IFS.
 The graphs of the two functions $L_0$ and $L_1$ are parallel lines.  When $\frac12 < a < 1$ the IFS ${\mathcal U}_a$ is overlapping. 
For the IFS ${\mathcal U}_{a} $ the
coding map will be denoted by $\pi_a$.  For a uniform IFS with mask point $p$,  the sections will be denoted by $\mu_{(a,p)}^+$ and  $\mu_{(a,p)}^-$, and  the associated dynamical systems by $([0,1], U_{(a,p)}^+)$ and $([0,1], U_{(a,p)}^-)$, where   $a < p<1-a$.  
\end{definition}

The following result concerning the uniform case follows readily from Parry \cite{parry}.  

\begin{theorem} \label{thmParry}
The topological entropy of the uniform dynamical systems $([0,1], U^{\pm}_{(a,q)})$ is equal to $-\ln(a)$.
\end{theorem}

\begin{lemma} \label{lemPi-omega}  
 Let $a \in (0,1)$ and $\omega \in \{0,1\}^{\infty}$.   For the IFS ${\mathcal U}_a$ we have 
$$\pi_a (\omega) = (1-a) \sum_{k=0}^{\infty} \omega_k \, a^k.$$  In particular, $\pi_a(\omega)$ is a continuous
function of $a$ in the interval $[0,1)$.
\end{lemma}

\proof
For the IFS ${\mathcal U}_a$ we have  $f_{i}(x)  = ax+i(1-a)$ for $i=0,1$.  Iterating
$$L_{\omega_0} \circ L_{\omega_1} \circ L_{\omega_{2}} \circ \cdots \circ L_{\omega_k}(x) = 
a^k x + (a^{k-1} \omega_{k-1} + \cdots + a\omega_1 + \omega_0)(1-a).$$
Therefore 
$$\pi_a (\omega) = \lim_{k \rightarrow \infty} L_{\omega|_k}(x)  = (1-a) \sum_{k=0}^{\infty} \omega_k \, a^k$$
for any $x\in [0,1]$.  Clearly the series converges for $0\leq a < 1$, and it is continuous inside the radius of 
convergence. 
\qed \m

\section{The Address Space} \label{secAddSpace}

The lexicographic order $\preceq$ on $\{0,1\}^{\infty}$ is the total order defined by 
 $\sigma \prec\omega$ if $\sigma\neq\omega$ and
$\sigma_{k}<\omega_{k}$ where $k$ is the least index such that $\sigma_{k}
\neq\omega_{k}$. For $\sigma,\omega\in \{0,1\}^{\infty}$ with $\sigma\preceq\omega$,
define the interval
$$\lbrack\sigma,\omega]    :=\{\zeta\in \{0,1\}^{\infty}:\sigma\preceq\zeta
\preceq\omega\},$$
and similarly for $(\sigma,\omega), \,
(\sigma,\omega]$, and $[\sigma,\omega)$.  We use the notation $\overline{0} = 000\cdots$ and $\overline{1} = 111\cdots$. 
Greek letters, other than coding map $\pi$ and section $\tau$, will denote strings; lower case Roman letters will denote
real numbers. Two itineraries play a special role. 

\begin{definition} \label{defMaskSeq}
For an overlapping masked IFS ${\mathcal F}$, the itineraries   $$\alpha _q := \tau^{-}_{q}(q) \, \qquad \qquad \text{and} \qquad \qquad \beta_q := \tau^{+}_{q}(q)$$ will be called the {\bf critical itineraries}.
\end{definition}

\begin{theorem} \label{thmAddSpace}  For an overlapping masked IFS ${\mathcal F}$ with mask point $q$, let $\overline{\Omega}_q = \Omega_q^+ \cup \Omega_q^-$.
\begin{enumerate}
\item  if $x,y \in [0,1]$ and $x > y$, then $(\tau_q^-)( x) \succ (\tau_q^+)( y)$;
\item the sections $\tau^+_q \, : \, [0,1]\rightarrow \Omega^+$ and $\tau^-_q \, : \, [0,1]\rightarrow \Omega^-$ are strictly increasing functions;
\item $\Omega^-_{q} = \{ \omega  \in \{0,1\}^{\infty} \, : \, S^n(\omega) \in [\overline{0},\alpha_q] \cup (\beta_q,\overline{1}] \quad$ for all $\quad n\geq 0 \}$;
\item $\Omega^+_{q} = \{ \omega  \in \{0,1\}^{\infty}  \, : \,  S^n(\omega) \in [\overline{0},\alpha_q) \cup [\beta_q,\overline{1}] \quad$ for all $\quad n\geq 0\}$;
\item $ \overline{\Omega}_{q} = \{ \omega  \in \{0,1\}^{\infty} 
 \, : \,  S^n(\omega) \in [\overline{0},\alpha_q] \cup [\beta_q,\overline{1}] \quad$ for all $\quad n\geq 0\}$.
\item  $\overline{\Omega}_q$ is the closure
of $\Omega_q^+$ and the closure of $\Omega_q^-$ in the metric space $\{0,1\}^{\infty}$.
\end{enumerate}
\end{theorem}

\proof  Since the mask is fixed, we suppress the index $q$ throughout the proof.  Also, when the superscript  $+$ or $-$ is omitted, we mean either one.  

Concerning statement 1, if $x > y$, then $(T^-)( x) > (T^+)( y)$ as long as
$x,y \leq q$ or $x,y \geq q$.  Hence $x > y$ implies that $\tau^-(x) \succeq \tau^+(y)$.  If $\tau^-(x) = \tau^+(y)$,
then $x = \pi(\tau^-(x)) = \pi(\tau^+(y)) = y$, a contradiction.  

Statement 2 follows directly from statement 1 since $\tau^+(x) \geq \tau^-(x)$ for all $x \in [0,1]$.

We next prove statement 3; the proof of statement 4 is omitted since it is done in essentially the same way.  To show that $\Omega^-_{q}$ is contained in $\{ \omega  \in \{0,1\}^{\infty} \, : \, S^n(\omega) \in [\overline{0},\alpha] \cup (\beta,\overline{1}] \quad \text{for all} \quad n\geq 0 \}$,  assume that $\omega \in \Omega^-$, and hence that $\omega = \tau^-(x)$ for some $x$.  If $\omega$ begins with a $0$, then $x \leq q$, which by the monotonicity of $\tau^-$ implies that $\omega = \tau^-(x) \preceq \tau^-(q) = \alpha$.  If $\omega$ begins with $1$, then $x > q$, which implies, using statement 1, that $\omega = \tau^-(x) \succ \tau^+(q) = \beta$. By shift invariance of $\Omega^-$,
the shift $S\omega \in \Omega^-$ and the same argument shows that $S \omega$ lies in the set $\{ \omega  \in \{0,1\}^{\infty} \, : \, S^n(\omega) \in [\overline{0},\alpha] \cup (\beta,\overline{1}] \quad \text{for all}\quad n\geq 0 \}$.

To prove  containment in the other direction in statement 3,  assume that $\omega = \omega _0 \omega_1 \omega_2 \cdots \in  \{ \omega  \in \{0,1\}^{\infty} \, : \, S^n(\omega) \in [\overline{0},\alpha] \cup (\beta,\overline{1}] \quad \text{for all}\quad n\geq 0 \}$.   By definition $\omega \in \Omega^-$ if $\omega$ lies in the image of $[0,1]$ under the section map.  By the definition of the section map, it is then sufficient to show that $M_{\omega|_k} \neq \emptyset$ for all $k$. We will show more, namely that $M_{(S^n \omega)|_k)} \neq \emptyset$ for all $k$ and all $n$.   This will be done by induction on $k$.   The statement is obviously true for $k = 0$.  Assuming it true for $k$, we will prove it for $k+1$.  Fix $n$ and let $j\sigma = (S^n  \omega)|_{k+1}$.   There are two cases, $j = 0$ and $j=1$. We will let $j=0$; the proof for $j=1$ is essentially the same.   By Lemma~\ref{lemM} it is sufficient to show that $ M_0 \cap f_0(M_{\sigma}) = M_{j \sigma} \neq \emptyset$.  Equivalently it must be shown that there is an $x\in M_{\sigma}$ such that $f_0(x) \leq q$.  By the induction hypothesis $M_{\sigma} \neq \emptyset$.  Since $0\sigma \preceq \alpha|_k$ and $\alpha_0=0$, also
$\sigma \preceq \widehat{\alpha} := \alpha_1\alpha_2\cdots \alpha_k$.  Since $\alpha \in \Omega^-$, we know that
$S^n  \alpha \in \Omega^-$, and hence $M_{\widehat{\alpha}} \neq \emptyset$, which implies that
there is a $y\in M_{0\widehat{\alpha}}$ such that $f_0(y) \leq q$.  But it follows easily from the definition of the 
partition $M^k$ that if $\sigma \preceq \widehat{\alpha}$, then the interval  $M_{\sigma}$ precedes (or is equal to)
the interval $M_{\widehat{\alpha}}$.  Therefore there is an $x\in M_{\sigma}$ such that $x \leq y$.  Since $f_0(y) \leq q$ and  $f_0$ is an increasing function, we arrive at the required $f_0(x) \leq q$.  

To prove statement 5, let $\Gamma = \{ \omega  \in \{0,1\}^{\infty} 
 \, : \,  S^n(\omega) \in [\overline{0},\alpha_q] \cup [\beta_q,\overline{1}] \;$ for all $\quad n\geq 0\}$.  Clearly
$\Omega^+ \subseteq \Gamma$ and $\Omega^- \subseteq \Gamma$.  Conversely $\Gamma \subseteq 
\Omega^+\cup \Omega^-$ unless there is a $\sigma \in \Gamma$ and integers $m$ and $n$ such that $S^n \alpha$ and $S^m = \beta$. Depending on whether $n > m$ or $m>n$, this implies that there is an integer $k$ such that
$S^k(\alpha) = \beta$ or $S^k(\beta) = \alpha$.  To show, by contradiction, that neither of  these equalities are possible, assume that $S^k(\alpha) = \beta$.  Since $\alpha \in \Omega^-$ and $\Omega^-$ is shift invariant,
also $\beta = S^k(\alpha) \in \Omega^-$.  But this contradicts the characterization of $\Omega^-$ given in statement 4. The equality  $S^k(\beta) = \alpha$ is likewise contradicted.

Statement 6 follows from statements 3, 4, and 5. 
\qed 

\begin{lemma} \label{lemContinuous}  For a masked overlapping IFS, the section $\tau_q^+ \, : \,[0,1] \rightarrow
\Omega^+_q$ is continuous at all points except those in the set $X^+ := \{x \, : \, S^n(\tau_q^+(x)) = \beta_q \; \text{for some} \; n\}$, and is continuous from the right everywhere.  Moreover, if $x\in X^+$ and $n$ is the least
integer such that  $S^n(\tau_q^+(x) ) = \beta_q$, then 
$$\lim_{y\rightarrow x^-} \tau^+_q(y) = \tau^+_q(x)|_n \, \alpha. $$
Likewise,  the section $\tau_q^- \, : \,[0,1] \rightarrow
\Omega^-_q$ is continuous at all points except those in the set $X^- := \{x \, : \, S^n(\tau_q^-(x)) = \alpha_q  \; \text{for some} \; n\}$, and is continuous from the left everywhere. Moreover, if $x\in X^-$ and $n$ is the least
integer such that  $S^n(\tau_q^+(x)) = \alpha_q $, then 
$$\lim_{y\rightarrow x^+} \tau^-_q(y) = \tau^-_q(x)|_n \, \beta. $$
\end{lemma}

\proof  To simplify notation, the subscript $q$ is omitted.  Consider the section $\tau^+$; the statement for $\tau^-$ is proved similarly.  The continuity at points not in $X^+$ follows directly from the continuity of $f_0$ and $f_1$ and the fact that $\tau^+$ can be viewed as an itinerary as described in Proposition~\ref{propItinerary}, likewise for the continuity from the right for points in $X^+$.  
From the definition of the dynamical system associated with the IFS, it is easy to verify that the following diagram commutes.
\begin{equation} \label{diagT}
\begin{array}
[c]{ccc}
[0,1] & \overset{T^{\pm}}{\rightarrow} & [0,1] \\
\tau^{\pm}\downarrow\text{\ \ \ \ } &  & \text{ \ \ \ }\downarrow\tau^{\pm}\\
\Omega_F & \underset{S}{\rightarrow} & \Omega_F
\end{array}
\end{equation}
By the commuting diagram above  $\tau^+((T^+)^n(x)) = S^n(\tau^+(x)) = \beta$, which implies that $(T^+)^n(x) = q$.  Since $n$ is the first such integer and if $y$ is sufficiently close to $x$, then $\tau^+(y) |_n = \tau^+(x) |_n$.
If $y<x$, then $(T^+)^n(y) < (T^+)^n(x) = q$.  Now
$$\begin{aligned} \lim_{y\rightarrow x^-} \tau^+(y) &= \tau^+(x)|_n \, \lim_{y\rightarrow x^-}  \tau^+((T^+)^n y) = \tau^+(x)|_n \, \lim_{y\rightarrow ((T^+)^n x)^-}  \tau^+(y) \\ &=
\tau^+(x)|_n \, \lim_{y\rightarrow q^-}  \tau^+(y) = \tau^+(x)|_n \, \lim_{y\rightarrow q^-}  \tau^-(y)  = 
 \tau^+(x)|_n \,  \alpha,
\end{aligned}$$
the second to last equality because, for any $m$ the first $m$ entries in the itineraries of $\tau^-(y)$ and $\tau^+(y)$ are equal if $y$ is sufficiently close to (and to the left of) $x$.  
\qed \m

Two dynamical systems $(\mathbb X, T)$ and $(\mathbb Y, S)$ are {\it topologically conjugate} if there exists a homeomorphism 
$\phi: \mathbb X \to \mathbb Y$ such that $ T =\phi^{-1} \circ S \circ \phi$.  Note that conditions 1 and 3 of Theorem~\ref{thmFT-TC} below   alone provide a necessary and sufficient condition for the fractal transformation  from one overlapping IFS to another to  be a homeomorphism.  The condition is simply that the critical itineraries of the associated dynamical systems be equal.  

\begin{theorem}  \label{thmFT-TC}
Given two overlapping masked IFSs $\mathcal F$ and $\mathcal G$  with respective mask points $q$ and $p$, sections  $\tau^{\pm}_F$ and $\tau^{\pm}_G$, dynamical systems $T^{\pm}_F$ and $T^{\pm}_G$, and address spaces $\Omega^{\pm}_F$ and  $\Omega^{\pm}_G$, the following statements are equivalent.
\begin{enumerate}
\item The fractal transformations  $\pi_G\circ \tau^{\pm}_F$  and $\pi_F\circ \tau^{\pm} _G $ are homeomorphisms. 
\item The address spaces are equal: $\Omega^{+}_F = \Omega^{+}_G$ and  $\Omega^{-}_F = \Omega^{-}_G$.
\item $\tau^+_F (q) = \tau^+_{G} (p)$ and  $\tau^-_F(q) = \tau^-_{G} (p)$. 
\item The dynamical systems $T^{+}_F$ and $T^{+}_G$ are topologically conjugate, as are $T^{-}_F$ and $T^{-}_G$. 
\end{enumerate}
\end{theorem}

\proof  To simplify notation we omit the superscript $\pm$.  We will show that $1 \Rightarrow 4  \Rightarrow 3\Leftrightarrow 2
 \Rightarrow 1$.  

$(1 \Rightarrow 4)\;$  Assume that  $h := \pi_G\circ \tau_F$  is a homeomorphism. 
Since $h$ is bijective, $\Omega_F = \Omega_G$.   From the commuting diagram~\ref{diagT} above and the fact that $\pi_G = \tau_G^{-1}$ on $\Omega_{G}=\Omega_F$, we have another commutative diagram for $G$.
\begin{equation*}
\begin{array}
[c]{ccc}
\Omega_{F} & \overset{S}{\rightarrow} & \Omega_{F} \\
\pi_G  \downarrow\text{\ \ \ \ } &  & \text{ \ \ \ }\downarrow \pi_G \\
   {[0,1]} & \underset{T_g} {\rightarrow} & [0,1]
\end{array}
\end{equation*}
Combining the two commutative diagrams we arrive at $T_G\circ h = h \circ T_F$ or $T_G = h_{FG}T_F h^{-1}$.  

$(4 \Rightarrow 3)\;$  Let topologically conjugate dynamical systems $T_F$ and $T_G$ be related 
by $T_G\circ h = h\circ T_F$, where $h$ is a homeomorphism.  If $q$ is the mask point of $\mathcal F$ and
$p$ is the mask point of $\mathcal G$, we claim that $p = h(q)$.  Otherwise, $h\circ T_F$ is 
discontinuous in some neighborhood of $q$ while $T_G\circ h$ is continuous in some neighborhood of $q$, a contradiction.
Now $T^n_F(q) \geq q$ if and only if $T^n_G(p) = T^n_G (h(q)) = h (T^n_F(q)) \geq h(q) = p$.  This implies 
statement (3).  

$(3 \Leftrightarrow 2)\;$  That  $(3 \Rightarrow 2)$ follows directly from statements 3 and 4 of Theorem~\ref{thmAddSpace}.  The same statements imply that the largest element of $\Omega^-$ that
starts with $0$ is $\alpha$, and the smallest element of $\Omega^+$ that starts with $1$ is $\beta$.  Therefore $(2 \Rightarrow 3)$.

$(2 \Rightarrow 1) \;$  To simplify notation we omit the subscript $q$.  
Assuming (2), we will show that $\pi_G \circ \tau^+_F$ is a homeomorphism.  Essentially the same proof shows that  $\pi_G \circ \tau^-_F$ is a homeomorphism.  Since $(2\Rightarrow 3)$ we know that the critical itineraries $\alpha$ and $\beta$ of
$\mathcal F$ are equal to the respective critical itineraries of $\mathcal G$, and moreover, for mask point $p$,
$$\pi_G (\alpha) = (\pi_G\circ \tau^-_G)(p)  = p = (\pi_G\circ \tau^+_G)(p)  = \pi_G (\beta).$$
Since it follows immediately from Definition~\ref{defFT} that  $\pi_G \circ \tau^+_F$ is a bijection, it suffices to show that it is continuous. (That the inverse in continuous is then a consequence of  Theorem~\ref{thmPi-Tau}.)   Because
$\pi_G$ is continuous,  Lemma~\ref{lemContinuous} implies  $\pi_G \circ \tau^+_F$ is continuous at all
points except perhaps those in the set  $X := \{x \, : \, S^n(\tau^+(x)) = \beta \; \text{for some} \; n\}$. 
Let $x \in X$.  Again by  Lemma~\ref{lemContinuous}, it suffices to prove that  $\pi_G \circ \tau^+_F$ is
continuous from the left.  But
$$\begin{aligned}  \lim_{y\rightarrow x^-} \pi_G (\tau^+_F(y)) &= \pi_G ( \lim_{y\rightarrow x^-} \tau^+_F(y)) =
\pi_G (\tau^+_F(x)|_n \alpha) = f_{\tau^+_F(x)|_n}(\pi_G\alpha) \\ &= f_{\tau^+_F(x)|_n}(\pi_G\beta) 
=\pi_G(\tau^+_F(x)|_n \beta ) =  \pi_G (\tau^+_F(x)) .\end{aligned}$$ 
\qed

\section{Entropy of an overlapping Dynamical System} \label{secEntropy}

Throughout this section, 
${\mathcal F}$ is an overlapping IFS with mask point $q$, critical itineraries $\alpha$ and $\beta$, and 
${\mathcal U}_a$ is a uniform IFS with coding map $\pi_a$.   In Lemmas~\ref{lemA} and \ref{lemC} we
assume that there exists an $a\in (0,1)$ such that  $\pi_a(\alpha) = \pi_a(\beta)$.  In this case
let 
\begin{equation} \label{eqC}  r(q) :=\min \, \left \{ a\in (0,1) \, : \,    \pi_a(\alpha) = \pi_a(\beta) \, \right \}.
\end{equation}
According to Lemma~\ref{lemPi-omega} 
\begin{equation} \label{eqD} r(q) = \min  \left \{ x\in (0,1) \, : \,   \sum_{n=0}^{\infty} \alpha_n x^n = \sum_{n=0}^{\infty} \beta_n x^n \, \right \}.\end{equation}
Note that the function $\sum_{n=0}^{\infty} (\beta_n -  \alpha_n) z^n$ is 
analytic   inside the unit disk in the complex plane, and hence can have at most finitely
many zeros within any closed disk of radius less than $1$. 
 In particular,
 $$\pi_r(\alpha) = \pi_r(\beta).$$ 

\begin{lemma} \label{lemA} (1)  Assume that there exists an $a\in (0,1)$ such that  $\pi_a(\alpha) = \pi_a(\beta)$ and let $r=r(q)$.  The map $\pi_a\, : \, \overline{\Omega}_q \rightarrow [0,1]$ is increasing for $0< a \leq r$ and
strictly increasing for $0< a < r$.    \m

(2) If there is no $a\in (0,1)$ such that $\pi_a(\alpha) = \pi_a(\beta)$, then the map $\pi_a \, : \, \overline{\Omega}_q \rightarrow [0,1]$ is strictly  increasing for all $a\in (0, 1)$.
\end{lemma}

\proof Since it is fixed throughout the proof, the subscript $q$ is omitted.
 Note that, for any IFS ${\mathcal U}_a$ with $a < \frac12$, it is easy to check, either because the
attactor is totally disconnected or directly from the power series, that if $\sigma \prec \omega$, then
$\pi_a(\sigma) < \pi_a(\omega)$.  

Concerning statement 1, let
\begin{equation} \label{eq-c}  s =\inf\{a \in(0,1)\, :\,  \pi_{a}(\sigma)=\pi_{a}(\omega) \;\; \text{for some} \;\;
 \sigma, \omega\in\overline{\Omega}, \, \sigma_0 \neq \omega_0  \}.
\end{equation}
Note that $s \leq r$ because $\alpha_0 = 0, \beta_0 = 1$ and $\alpha, \beta \in \overline{\Omega}$.  
Using the  continuity of $\pi_{a}(\sigma)$ in $a$  (see
Lemma~\ref{lemPi-omega}) and $\sigma$  (see the comments following Definition~\ref{defCode}), and the
compactness of $\overline{\Omega}$,  it follows that there exist $\sigma,\omega
\in\overline{\Omega}$ such that $\pi_{s}(\sigma)=\pi_{s}(\omega)$.
We claim that $r = s$.  Assume, by way of contradiction, that $s < r$.  If we assume, without loss
of generality that $\sigma_0 = 0$ and $\omega_0=1$, then
\begin{equation} \label{eq-a}
\pi_{\frac13}(\sigma)\leq\pi_{\frac13}(\alpha)<\pi_{\frac13}(\beta)\preceq \pi_{\frac13}(\omega)
\end{equation}
because by Theorem~\ref{thmAddSpace} we have $\sigma \preceq \alpha$ and $\beta \preceq \omega$ and, as mentioned above, $\pi_{\frac13}$ is order preserving.  Consider $\pi_a(\sigma), \pi_a(\alpha),
\pi_a(\beta), \pi_a(\omega)$ as functions of $a \in [1/3,r]$. (It is helpful to visualize the graphs of these
these four functions.)   Since $s<r$, we have
\begin{equation} \label{eq-b} \begin{aligned} \pi_a(\alpha) < \pi_a(\beta) &\quad \text{for} \quad   \frac13 \leq a \leq s \\
\pi_a(\sigma) < \pi_a(\omega)  &\quad \text{for} \quad  \frac13 \leq a < s \\
\pi_{s}(\sigma)=\pi_{s}(\omega) &
\end{aligned}
\end{equation}
By the continuity of $\pi_a$ with respect to $a$ and the intermediate value theorem,
the formulas~\ref{eq-a} and \ref{eq-b} imply that either there is a $t\in (\frac13, s)$ such that
$\pi_t(\sigma) = \pi_t(\alpha)$  with $\sigma \neq \alpha$ or there is a  $t\in (\frac13, s)$ such that $\pi_t(\omega) = \pi_t(\beta)$ with $\omega \neq \beta$.  Since the proof is essentially the same in either
case, assume that $\pi_t(\sigma) = \pi_t(\alpha)$  with $\sigma \neq \alpha$.  Since $t < s$, this
would contradict the minimality of $s$ (in eqaution~(\ref{eq-c})) if $\sigma_0 = 0$ and $\alpha_0=1$.
This is not the case, however, because $\alpha_0 = 0$.  In order to get the contradiction, we define two related strings
$\sigma'$ and $\alpha'$ such that  $\pi_t(\sigma') = \pi_t(\alpha')$ and  $\sigma'_0 = 0$ and $\alpha'_0=1$.
To do this,  let $k$ be the least integer such that $(S^{k}\sigma)_{0}\neq(S^{k} \alpha)_{0}$ and let $\sigma' = S^{k}\sigma$
and $\alpha' = S^{k}\alpha$, which forces $\sigma_{0} \neq \omega'_0$.  We are now done
because $ L_{\sigma|_k} (\pi_t(\sigma')) = \pi_t(\sigma) =  \pi_t(\omega) = L_{\omega|_k} (\pi_t(\omega')) =
 L_{\sigma|_k} (\pi_t(\omega'))$ implies, because $ L_{\sigma|_k}$ is invertible, that  $\pi_t(\sigma') =  \pi_t(\omega')$.
The shift invariance of $\overline{\Omega}$ guarantees that $\sigma', \omega' \in \overline{\Omega}$.
Therefore $s=r$.  

To conclude the proof of statement 1 of the lemma, assume that   $a \prec r$, $\sigma,\omega\in\overline{\Omega}$, and $\sigma\prec\omega$.   If $\sigma_0 = 0$ and $\omega_0 = 1$, then 
$\pi_{a}(\sigma)\neq\pi_{a}(\omega)$ by what was proved in the paragraph above.  Since
$\pi_{\frac13}(\sigma) < \pi_{\frac13}(\omega)$, it would follow that $\pi_{a}(\sigma) < \pi_{a}(\omega)$;
otherwise the crossing graphs would contradict $s=r$.  Even if  $\sigma_0 = \omega_0$, we claim that
$\pi_{a}(\sigma)\neq\pi_{a}(\omega)$.  Assume otherwise, that $\pi_{a}(\sigma) =\pi_{a}(\omega)$,
then by letting $\sigma'$ and $\omega'$ be shifts of $\sigma$ and $\omega$, respectively, exactly as was
done in the paragraph above, we get $\pi_{a}(\sigma') =\pi_{a}(\omega')$ with $\sigma'_0 \neq \omega'_0$,
which contradicts $s=r$. 

 In the case $a = r$ and $\sigma\prec\omega,$ clearly $\pi_{a}(\sigma) > \pi_{a}(\omega)$  
could contradict the continuity of $\pi_a$ at $a=r$; therefore  $\pi_{a}(\sigma) \leq \pi_{a}(\omega)$ .

Lastly consider statement 2, i.e. the case $\pi_{a}(\alpha) \neq \pi_a(\beta)$ for all $a \in (0,1)$.  Essentially the same proof as above shows that $s=1$ and consequently that if $\sigma \prec \omega$ then  $\pi_{a}(\sigma)  < 
\pi_{a}(\omega)$ for all $a < s=1$.
\qed 

\begin{lemma} \label{lemC} Assume that there exists an $a\in (0,1)$ such that  $\pi_a(\alpha) = \pi_a(\beta)$ and let $r=r(q)$.  For any integer $n>0$, if $S^n(\beta) \prec \alpha$, then  $\pi_r(S^n(\beta)) < \pi_r(\alpha)$.  Similarly
if  $S^n(\alpha) \succ \beta$, then $\pi_r(S^n(\alpha)) > \pi_r(\beta)$.
\end{lemma}

\proof 
The following are readily verifiable facts about the partitions of $[0,1]$ that are part of Definition~\ref{defSM} of the sections associated with the masks $M^+_q$ and $M^-_q$.   Denote the $k^{th}$ partitions by $(M^k)^+$ and $(M^k)^-$.  
\begin{enumerate}
\item The sets in partitions  $(M^k)^+$ (except the last)  and $(M^k)^-$  (except the first) are half open intervals of the form $[\cdot , \cdot)$ and $( \cdot , \cdot]$, respectively.
\item The endpoints of the intervals  in $(M^k)^+$ have the same endpoints as the intervals in $(M^K)^-$.  Denote the
set of open intervals by $M^k$.  
\item Given any interval $I$  in $M^k$, the first $k$ elements in the address (either $+$ or $-$ address) of any two points in $I$  are equal. 
\item If $(x,y)$ is an interval in $M^k$ whose elements have address beginning with $\theta \theta_1 \theta_2$,
where $\theta$ has length $k-2$ and $\theta_1, \theta_2 \in \{0,1\}$,  then the address $\tau^+(x)$ of $x$ is
$$\begin{cases}  \theta 0 \beta  \quad \text{if} \quad \theta_1\theta_2 = 01 
\\ \theta \beta \quad \;\text{if} \quad \theta_1\theta_2 = 10, 
\end{cases}$$
and the address $\tau^-(y)$ of $y$ is
$$\begin{cases} \theta \alpha \quad \; \text{if} \quad \theta_1\theta_2 = 01 
\\ \theta 1\alpha \quad \text{if} \quad \theta_1\theta_2 = 10. 
\end{cases}$$
\end{enumerate}

 We will prove that  $S^n(\beta) \prec \alpha$ implies  $\pi_r(S^n(\beta)) < \pi_r(\alpha)$.  That
 $S^n(\alpha) \succ \beta$ implies $\pi_r(S^n(\alpha)) < \pi_r(\beta)$ has essentially the same proof.
Assume that $S^n(\beta) \prec \alpha$.  There exists a $k$ (sufficiently large) and three open intervals $I_1 = (x_1,y_1), \, I_2=(x_2,y_2),\,  I_3 = (x_3,y_3)  \in M^k$ with the following properties: \m

5.  $y_1 \leq x_2  < y_2 \leq x_3 $, \m

6.  $\tau(z)|_k = \alpha|_k$ for all $z\in I_3$, \m

7.  $y_3 = q$, \m

8.  $\tau(z)|_k = S^n(\beta)|_k$ for all $z\in I_1$,  \m 

9.  either the last two elements $\tau(z)|_k$  are $01$ for all  $z \in I_2$, or the last two elements $\tau(z)|_k$  are $10$ for all  $z \in I_2$, and \m

10.  $S^n(\beta) = \tau^+(z_0)$ for some $z_0\in [x_1,y_1)$. \m

\noindent  The existence of the intermediate interval $I_2$ follows from the facts that the right endpoint $y_3$ of $I_3$ is fixed at $q$ (statement 7) and that the lengths of the intervals of $M^k$ tends to $0$ as $k\rightarrow \infty$.  
If statement 9 were false, then there would  exist a $k$, an interval $I \in M^k$, and a finite string $\theta$ such that $\tau(z) = \theta\overline{0}$ or $\tau(z) = \theta\overline{1}$ for all $z\in I$, which is impossible  
 (again because the lengths of the intervals of $M^k$ tends to $0$ as $k\rightarrow \infty$).

By statement 2 of Theorem~\ref{thmAddSpace} and by propertiy 4 above, if the last two elements of the finite $I_2$-address is $01$, then
\begin{equation} \label{eqAA}
S^n(\beta) =   \tau^+(z_0)  \preceq \tau^+(x_2) = \theta 0\beta \qquad \text{and} \qquad \theta \alpha = \tau^-(y_2) \preceq \tau^{-}(y_3) = \alpha.\end{equation} 
If the last two elements of the finite $I_2$-address is $10$, then 
\begin{equation} \label{eqBB} S^n(\beta) =   \tau^+(z_0)  \preceq \tau^+(x_2)  = \theta \beta \qquad \text{and} \qquad \theta 1 \alpha= \tau^-(y) \preceq \tau^{-}(y_3) = \alpha.\end{equation}

\noindent Consider the first case above; the proof for the second case is essentially the same. From the inequalities
above and by Lemma~\ref{lemA} (since $\theta 0 \beta$ and $\theta \alpha$ lie in $\overline{\Omega}$), we have $\pi_r(S^n(\beta))  \leq \pi_r(\theta 0\beta)$ and $\pi_r(\theta \alpha) \leq \pi_r(\alpha)$. The proof is complete if 
$ \pi_r(\theta 0\beta) < \pi_r(\theta \alpha) $.  But using Lemma~\ref{lemPi-omega}
$$\begin{aligned}
\pi_r(\theta\alpha) - \pi_r(\theta 0\beta) &= r^k(\pi_r(\alpha)-\pi_r(0\beta) )= r^k(1-r)(\pi_r(\alpha) - r\pi_r(\beta)) \\
& =r^k(1-r)(\pi_r(\alpha) - r\pi_r(\alpha) )= r^k(1-r)\pi_r(\alpha) > 0, \end{aligned}.$$
\qed

\begin{lemma} \label{lemD}  There exists an $a\in (0,1)$ such that $\pi_a(\alpha) = \pi_a(\beta)$.
\end{lemma}

\proof  Assume, by way of contradiction, that $\pi_a(\alpha) < \pi_a(\beta)$ for all $a \in (0,1)$.  Let $a$
be arbitrary in the interval $(0,1)$.  By statement 2 of Lemma~\ref{lemA},  the map $\pi_a$ is strictly increasing  on $\overline{\Omega}_q$.  Let $p = \pi_a(\beta)$.   For $\omega \in \Omega^+_q$ we claim that 
$$U^+_{(a,p)}(\pi_a\omega) = \pi_a(S\omega),$$ where  $U^+_{(a,p)}$ is the uniform dynamical system.  
This would  imply that the address space $\Omega^+_q$ is an invariant subset of the dynamical system $U^+_{(a,p)}$.  This, in turn, would  imply that the entropy of the overlapping dynamical system $T^+_q$ is less than or equal 
to the entropy of the uniform dynamical system $U^+_{(a,p)}$, which, according to Theorem~\ref{thmParry}, equals $-\ln a$.  Since this is true 
for all $a\in (0,1)$, the entropy of  $T^+_q$ must be $0$, which is not possible for a dynamical system where
where the two continuous branches are expansive. 

To prove the claim,  let $U = U^+_{(a,p)}$.  First note, from Lemma~\ref{lemA}, that if
$\pi_a(\omega) < p = \pi_a(\beta)$ then $\omega \prec \beta$, and hence $\omega_0=0$.
Likewise if $\pi_a(\omega) \geq p = \pi_a(\beta)$ then $\omega \succeq \beta$, and hence $\omega_0=1$.
Therefore if $\pi_a(\omega) < p$ then
$$U(\pi_a(\omega)) = U( (1-a) \sum_{n=0}^{\infty} \omega_n a^n)  = 
(1-a) \sum_{n=0}^{\infty} \omega_{n+1} a^n  + \frac{\omega_{0}}{a} = 
(1-a) \sum_{n=0}^{\infty} \omega_{n+1} a^n  =  \pi_a(S \omega),$$ and if 
$\pi_a(\omega) \geq p$, then
$$\begin{aligned} U(\pi_a(\omega)) &= U ( (1-a) \sum_{n=0}^{\infty} \omega a^n )  = 
(1-a) \sum_{n=0}^{\infty} \omega_{n+1} a^n  + \frac{\omega_{0}}{a} - \frac{1}{a} \\
&=  (1-a) \sum_{n=0}^{\infty} \omega_{n+1} a^n  =  \pi_a (S\omega).\end{aligned}$$
\qed

\begin{lemma} \label{lemB}  Let $r = r(q)$ and 
 $p = \pi_r(\alpha)= \pi_r(\beta)$.  If the uniform IFS ${\mathcal U}_r$  with coding map $\pi_r$, has mask point $p$ and sections $\mu^{+}_{(r,p)}$  and $\mu^{-}_{(r,p)}$, then $\mu^{-}_{r,p}(p) = \alpha$ and   $\mu^{+}_{r,p}(p) = \beta$.
\end{lemma}

\proof We will prove that  $\mu^{+}_{(r,p)}(p) = \beta$; the proof that $\mu^{-}_{(r,p)}(p) = \alpha$  is 
essentially the same.   Let  $U := U^+_{(r,p)}$  be the dynamical system associated with
the uniform IFS and let $\omega := \mu^{+}_{(r,p)}(p)$.  For all $n\geq 0$, we will prove the following by induction on $n$:
\begin{enumerate}
\item $\omega_{n} = \beta_{n}\quad$ and,
\item $U^n(p) = \pi_r(S^n \beta)$.
\end{enumerate}
Since both $\beta$ and $\omega$ begin with a $1$, both statements are true for $n=0$.  
Assuming  the two statements true for $n-1$, we will prove that they are true for $n$.   

Starting with statement 2: 
$$\begin{aligned} U^{n}(p) &= U(U^{n-1}(p) ) = U(\pi_r(S^{n-1} \beta) )
 =  U ((1-r) \sum_{k=0}^{\infty} \beta_{n-1+k} r^k )  \\  &=
(1-r) \sum_{k=0}^{\infty} \beta_{n+k} r^k  = \pi_r (S^n \beta). \end{aligned}$$   The second to
last equality above comes from the following direct calculation:  if $\omega_{n-1} = 0$, then by the induction
hypothesis $\beta_{n-1}=0$ and
$$U( (1-r) \sum_{k=0}^{\infty} \beta_{n-1+k} r^k )  = 
(1-r) \sum_{k=0}^{\infty} \beta_{n+k} r^k  + \frac{\beta_{n-1}}{r} = 
(1-r) \sum_{k=0}^{\infty} \beta_{n+k} r^k  =  \pi_r (S^n \beta),$$ and if 
$\omega_{n-1} = 1$, then $\beta_{n-1} =1$ and 
$$U ( (1-r) \sum_{k=0}^{\infty} \beta_{n-1+k} r^k )  = 
(1-r) \sum_{k=0}^{\infty} \beta_{n+k} r^k  + \frac{\beta_{n-1}}{r} - \frac{1}{r} = 
(1-r) \sum_{k=0}^{\infty} \beta_{n+k} r^k  =  \pi_r (S^n \beta).$$

Concerning statement 1, 
if $\beta_n = 0$, then by statement 4 of Theorem~\ref{thmAddSpace} we have $S^n \beta \prec \alpha$.  Therefore, by Lemma~\ref{lemC} and statement 1 which we have just proved, we have  $U^n(p)   =  \pi_r(S^n \beta) < \pi_r (\alpha) = p$.   By the definition of the itinerary of $p$ this implies that  $\omega_n = 0$, and hence $\omega_n = \beta_n$.   If, on the other hand, $\beta_n = 1$, then by statement 3 of Theorem~\ref{thmAddSpace} we have $S^n \beta \succeq \beta$, and therefore  $U^n(p)   =  \pi_r(S^n \beta) \geq \pi_r (\beta) = p$.   Again by definition of the itinerary of $p$, we have $\omega_n = 1$ and hence $\omega_n = \beta_n$.   
\qed 

\begin{theorem}  \label{thmMain}  Let $([0,1],T)$ be any overlapping dynamical system with mask point $q$, critical itineraries $\alpha$ and $\beta$, and $r(q)$ as defined in equations~(\ref{eqC}) or (\ref{eqD}).
\begin{enumerate}
\item The  dynamical system $([0,1], T)$ is topologically conjugate to the uniform dynamical system $([0,1), U_{r,p})$, where $r = r(q)$ and  $p= (1-r) \sum_{n=0}^{\infty} \alpha_n r^n$.  
\item The entropy of dynamical system $([0,1], T)$ is $-\ln r$, where $r$ the smallest solution $x\in [0,1]$ to the equation
$$\sum_{n=0}^{\infty} \alpha_n x^n = \sum_{n=0}^{\infty} \beta_n x^n.$$
\end{enumerate}
\end{theorem}

\proof Statement 1 follows immediately from Lemma~\ref{lemB} and from the equivalence of statements 3 and 4 
of Theorem~\ref{thmFT-TC}.  Statement 2 then follows immediately from Theorem~\ref{thmParry} and the fact that
two topologically conjugate dynamical systems have the same entropy.
\qed \B

{\bf Acknowledgement.}  We thank Jed Keesling, Nina Snigireva, and Tony Samuels for many helpful discussions in connection with this work.  We thank Konstantin Igudesman for useful comments on an early version of this work. \B

\end{document}